\documentclass[12pt,a4paper]{article}
\usepackage{amsmath,amssymb,amsthm}
\usepackage{amssymb,latexsym, bbm,comment}
\usepackage{graphicx}
\usepackage{makeidx}

\newtheorem{theorem}{Theorem}[section]

\newtheorem{prop}{Proposition}[section]

\newcommand{\nN}{n \in \mathbb{N}}
\newcommand{\N}{\mathbb{N}}
\newcommand{\C}{\mathbb{C}}
\newcommand{\E}{\mathbb{E}}

\newcommand{\ds}{\displaystyle}

\newcommand{\tr}{\rm{tr}}

\newcommand{\bean}{\begin{eqnarray*}}
\newcommand{\eean}{\end{eqnarray*}}

\newcommand{\Z}{\mathbb{Z}}

\newcommand{\G}{\widehat{G}}
\newcommand{\GF}{\widehat{G_{F}}}

\makeindex

\begin{document}

\date{}

\title{Stationary Random Fields on the Unitary Dual of a Compact Group }

\author{ David Applebaum, \\ Department of Probability and Statistics,\\ University of
Sheffield,\\ Hicks Building, Hounsfield Road,\\ Sheffield,
England, S3 7RH\\ ~~~~~~~\\e-mail: D.Applebaum@sheffield.ac.uk }

\maketitle

\begin{abstract} We generalise the notion of wide-sense stationarity from sequences of complex-valued random variables indexed by the integers, to fields of random variables that are labelled by elements of the unitary dual of a compact group. The covariance is positive definite, and so it is the Fourier transform of a finite central measure (the spectral measure of the field) on the group. Analogues of the Cramer and Kolmogorov theorems are extended to this framework. White noise makes sense in this context and so, for some classes of group, we can construct time series and investigate their stationarity. Finally we indicate how these ideas fit into the general theory of stationary random fields on hypergroups.
\end{abstract}

\maketitle

\section{Introduction}

There are many important classes of stochastic process that have been systematically developed, both because of their mathematical vitality, and their importance for applications. These include, for example,  Markov chains, branching processes, and diffusion processes. The emphasis in this paper is on discrete-time (wide-sense) stationary, complex-valued processes $(X_{n}, n \in \Z)$, so that $\E(|X_{n}|^{2}) < \infty$ and
\begin{equation} \label{statdeford}
\E(X_{m} \overline{X_{n}}) = \E(X_{m-n}\overline{X_{0}}),
\end{equation}
for all $m,n \in \Z$. These may be used to model fluctuations from some fixed background signal. Stationarity is a vital ingredient in the theory of time series (see e.g. \cite{BD}) which has a wide range of applications, including economics and climate science.

A stationary process is characterised by its covariance function $C(n) = \E(X_{n}\overline{X_{0}})$, which is positive-definite, and so by the Herglotz theorem, there is a finite measure $\mu$ on the torus $\mathbb{T}$, known as  the spectral measure of the process, for which
$$ C(n) = \int_{\mathbb{T}}e^{-2\pi i n \theta} \mu(d \theta),$$
for all $\nN$.

If we are interested in describing the interaction of chance with symmetry, then it is natural to consider stationary random fields on a group $G$, i.e. mappings $X:G \rightarrow L^{2}(\Omega, \C)$ for which\footnote{We only write down the left-invariant case here, but of course right-invariance is equally valid.}
$$ \E(X(hg_{1})\overline{X(hg_{2})}) = \E(X(g_{1})\overline{X(g_{2})}),$$
for all $g_{1}, g_{2}, h \in G$. The study of these, and related objects on homogeneous spaces, seems to have begun with work by A.M. Yaglom in the late 1950s (see e.g. \cite{Yag}); recently there have been monograph treatments and new applications to e.g. earthquake modelling and the study of the cosmic background radiation left over from the Big Bang (\cite{Maly, MaPe}).

In this paper, we suggest that, although replacing $\Z$ as the index of a stationary field by a group $G$ is mathematically highly productive, it may not be the most natural generalisation. As was discussed above, the spectral measure of a stationary process is defined on the torus $\mathbb{T}$; this is the simplest compact group, and its dual group is $\Z$. We propose that $\mathbb{T}$ should be replaced by a general (and so, not necessarily abelian) compact group, so that the role of $\Z$ is now played by the unitary dual $\G$ of $G$. Note that $\G$ is not itself a group if $G$ fails to be abelian.

In section 2 of this paper we generalise the definition (\ref{statdeford}) to random fields over $\G$. Indeed we say that a field $(Y_{\pi}, \pi \in \G)$ is stationary if
$$ \E(Y_{\pi_{1}}\overline{Y_{\pi_{2}}}) = \E(Y_{\pi_{1} \otimes \pi_{2}^{*}}\overline{Y_{\epsilon}}),$$
for all $\pi_{1}, \pi_{2} \in \G$, where $\pi^{*}$ is the irreducible representation that is conjugate to $\pi$, and $\epsilon$ is the trivial representation. Some justification as to why this is a sensible generalisation of (\ref{statdeford}) will be provided. We also define the covariance function and show that it is the Fourier transform (in the group-theoretic sense) of a finite central measure on $G$, which we call the spectral measure of the field.  We establish a Cramer-type representation of stationary fields as stochastic integrals with respect to orthogonally scattered random measures on $G$, and we prove a theorem of Kolmogorov-type to the effect that every positive-definite function on $\G$ is the covariance of a stationary random field on $\G$.

We have already pointed out that $\G$ is not in general a group, but it is a hypergroup (\cite{BH}) and we discuss this in section 4. There is an existing literature on stationary random fields on hypergroups \cite{Hesurv, LaLe, Leit} which this current work complements. We make some observations:
\begin{enumerate}
\item The definition of stationarity for general hypergroups is quite non-intuitive. But in our case, the parallel with the classical case is very direct.
\item The duality between the hypergroup $\G$ and the group $G$ is manifest in the relationship between the stationary field and its spectral measure. There is a rich structure here that merits further investigation, and which could lead to new examples of the important class of central measures on compact groups.
\item The key process of ``white noise '' may not exist in general hypergroups. But it always does in our case. This means that, at least for some classes of compact groups, we may develop a theory of time series on their unitary duals, and investigate stationarity. Some examples for the case of the dual of $SU(2)$ are considered in section 3 of this paper.
\end{enumerate}

{\it Notation}. If $A$ is  a complex-valued matrix, then tr$(A)$ is its trace (i.e. the sum across the leading diagonal). If $U$ is a topological space, then ${\cal B}(U)$ is the Borel $\sigma$-algebra of $U$ (i.e. the smallest $\sigma$-algebra containing all open sets). Haar integrals of suitable functions $f$ on a compact group $G$ are written $\int_{G}f(\sigma)d\sigma$.

\section{Definition and Main Results}

Let $G$ be a compact (second countable, Hausdorff) topological group, $\G$ be its unitary dual, comprising equivalence classes of irreducible unitary representations of $G$, and $\GF$ be the set of equivalence classes of finite-dimensional unitary representations of $G$ (each with respect to unitary conjugation). Since $G$ is compact, $\G$ is countable and $\G \subset \GF$.\footnote{We refer to a standard text, such as \cite{BroD}, for all facts about compact groups quoted herein. See also the account for probabilists in \cite{App}.} We denote the trivial representation of $G$ by $\epsilon \in \G$.  The {\it character} $\chi_{\pi}$ of $\pi \in \GF$ is defined by
$$ \chi_{\pi}(g) = \tr(\pi(g))$$ for each $g \in G$, and it is consequence of the celebrated Peter-Weyl theorem that $\{\chi_{\pi}, \pi \in \G \}$ is a complete orthonormal basis in the complex Hilbert space $L^{2}_{c}(G)$ of all central (i.e. conjugate invariant) square-integrable (with respect to normalised Haar measure) functions on $G$. It follows that we may decompose each $\pi \in \GF$ as \begin{equation} \label{multp} \pi = \bigoplus_{\pi^{\prime} \in \G}M(\pi, \pi^{\prime})\pi^{\prime}, \end{equation} where \begin{equation} \label{mult}M(\pi, \pi^{\prime}) = \int_{G}\chi_{\pi}(g^{-1})\chi_{\pi^{\prime}}(g)dg,\end{equation} is the {\it multiplicity} of $\pi^{\prime}$ is $\pi$. Of course $M(\pi, \pi^{\prime}) \in \Z_{+}$ vanishes for all but finitely many $\pi^{\prime} \in \G$.
The conjugate representation associated to $\pi \in \GF$ is denoted $\pi^{*}$ and the tensor product of the representations $\pi_{1}$ and $\pi_{2}$ is $\pi_{1} \otimes \pi_{2}$. Note that for all $g \in G$,
\begin{equation} \label{charp}
\chi_{\pi^{*}}(g) = \overline{\chi_{\pi}(g)}~~~,~~~ \chi_{\pi_{1} \otimes \pi_{2}}(g) = \chi_{\pi_{1}}(g)\chi_{\pi_{2}}(g),
\end{equation}
for each $g \in G$.

\begin{prop} \label{useful}
For all $\pi_{1}, \pi_{2} \in \G$,
$$ M(\epsilon, \pi_{1} \otimes \pi_{2}^{*}) = \delta_{\pi_{1}, \pi_{2}}.$$
\end{prop}

{\it Proof.} Using (\ref{mult}), (\ref{charp}), and orthonormality of characters, we have
\bean M(\epsilon, \pi_{1} \otimes \pi_{2}^{*})  & = & \int_{G} \chi_{\pi_{1}}(g)\overline{\chi_{\pi_{2}}}(g)dg\\
& = & \delta_{\pi_{1}, \pi_{2}}.~~~~~~~~~~~~~~~~~~~~~~ \Box  \eean

Let $(\Omega, {\cal F}, P)$ be a probability space. A mapping $Y: \GF \rightarrow L^{2}(\Omega, {\cal F}, P;\C)$ is said to be a {\it decomposable random field} on $\GF$ if it satisfies
$$ Y_{\pi} = \sum_{\pi^{\prime} \in \G}M(\pi, \pi^{\prime})Y_{\pi^{\prime}},$$
with respect to the decomposition (\ref{multp}). Clearly such a field is uniquely determined by its values on $\G$. We say that such a field is  {\it (wide-sense) stationary} if
\begin{equation} \label{statdef}
\E(Y_{\pi_{1}}\overline{Y_{\pi_{2}}}) = \E(Y_{\pi_{1} \otimes \pi_{2}^{*}}\overline{Y_{\epsilon}}),
\end{equation}
for all $\pi_{1}, \pi_{2} \in \G$. The motivation for the definition (\ref{statdef}) comes from the well-known case $G = \mathbb{T} = [0, 2\pi), \G = \Z$. In that case the irreducible representation corresponding to $\pi_{1} = n$ is uniquely determined by the character $\theta \rightarrow e^{i n\theta}$, and the character associated to $\pi_{1} \otimes \pi_{2}^{*}$, where $\pi_{2} = m$, is precisely  $\theta \rightarrow e^{i (n-m)\theta}$.
\vspace{5pt}

{\bf Remark}. Clearly if the random field is stationary, then $\E(|Y_{\pi}|) < \infty$, for all $\pi \in \G$. It may seem strange to some readers that we do not impose some additional stationarity condition on the means, i.e. that the quantity $\E(Y_{\pi})$ does not depend on $\pi \in \G$, or even that the field is {\it centred}, in that $\E(Y_{\pi}) = 0$, for all $\pi \in \G$. Here we follow Doob \cite{Doob}, p.95 who, in the classical case $G = \mathbb{T}, \G = \Z$ wrote, ``Usually the added condition that $\E(X_{s})$ does not depend on $s$ is imposed. This condition is unnatural mathematically, and has nothing to do with the essential properties of interest in these processes, and we shall therefore not impose it.''

\vspace{5pt}

If $Y$ is a stationary random field on $\G$, we define its {\it covariance function} $C_{Y}: \GF \rightarrow \C$ by
$$ C_{Y}(\pi) = \E(Y_{\pi}\overline{Y_{\epsilon}}),$$
for all $\pi \in \GF$, and we note that it is a decomposable mapping on $\GF$ in that
$$ C_{Y}(\pi) = \sum_{\pi^{\prime} \in \G}M(\pi, \pi^{\prime})C_{Y}(\pi^{\prime}),$$
with respect to (\ref{multp}).

We recall from \cite{Hey} that $\Phi:\G \rightarrow \C$ is {\it positive definite} if for all $N \in \N, \pi_{1}, \ldots, \pi_{N} \in \G$ and $c_{1}, \ldots, c_{N} \in \C$,
$$ \sum_{m, n = 1}^{N}c_{m}\overline{c_{n}}\sum_{\pi \in \G}M(\pi, \pi_{m} \otimes \pi_{n}^{*})\Phi(\pi) \geq 0.$$
If $\Phi$ extends to a mapping $\GF \rightarrow \C$ that is decomposable, then we have the equivalent condition
\begin{equation} \label{pd}
 \sum_{m, n = 1}^{N}c_{m}\overline{c_{n}}\Phi(\pi_{m} \otimes \pi_{n}^{*}) \geq 0.
\end{equation}

\begin{prop} \label{covpd}
If $Y$ is a stationary random field on $\G$, then its covariance function $C_{Y}$ is positive definite.
\end{prop}

{\it Proof.} Using (\ref{pd}) and (\ref{statdef}), we find that
\bean \sum_{m, n = 1}^{N}c_{m}\overline{c_{n}}C_{Y}(\pi_{m} \otimes \pi_{n}^{*}) & = & \sum_{m, n = 1}^{N}c_{m}\overline{c_{n}}\E(Y_{\pi_{m} \otimes \pi_{n}^{*}}\overline{Y_{\epsilon}})\\
& = & \sum_{m, n = 1}^{N}c_{m}\overline{c_{n}}\E(Y_{\pi_{m}}\overline{Y_{\pi_{n}}})\\
& = & \E\left(\left|\sum_{n=1}^{N}c_{n}Y_{\pi_{n}}\right|^{2}\right) \geq 0.~~~~~~\Box \eean

It follows from Proposition \ref{covpd} and the Bochner theorem (Theorem 5.5 in \cite{Hey}) that there exists a finite Radon central measure $\mu_{Y}$ defined on the Borel $\sigma$-algebra of $G$ for which
\begin{equation} \label{spec}
C_{Y}(\pi) = \int_{G} \chi_{\pi}(g) \mu_{Y}(dg)
\end{equation}
for all $\pi \in \G$. We call $\mu_{Y}$ the {\it spectral measure} of the random field $Y$.

As an example, consider the {\it white noise} $Z: \GF \rightarrow L^{2}(\Omega, {\cal F}, P)$ which is defined to be a decomposable random field which is {\it uncorrelated} in that
$$ \E(Z_{\pi}\overline{Z_{\pi^{\prime}}}) = \delta_{\pi, \pi^{\prime}},$$
for all $\pi, \pi^{\prime} \in \G$. It follows from Proposition \ref{useful} that $Z$ is stationary and the spectral measure is easily seen to be (normalised) Haar measure on $G$.

The next result gives a {\it Cramer representation} for the field.

\begin{theorem} If $(Y_{\pi}, \pi \in \G)$ is a stationary random field, then there exists an orthogonally scattered random field $\Gamma_{Y}$ on $G$ so that for all $\pi \in \G$,
\begin{equation} \label{cram1}
              Y_{\pi} = \int_{G}\chi_{\pi}(g)\Gamma_{Y}(dg)~\mbox{a.s.}.
\end{equation}
Furthermore, $\E(|\Gamma_{Y}(A)|^{2} = \mu_{Y}(A)$ for all $A \in {\cal B}(G)$, and $\Gamma_{Y}$ is a.s. central in that for each $g \in G$,
$$ P(\Gamma_{Y}(gAg^{-1}) = \Gamma_{Y}(A)) = 1.$$

\end{theorem}

{\it Proof.} This is along standard lines. We sketch the details following the argument given in \cite{Lamp} pp. 46-7 for the classical case.

Let ${\cal M}$ be the closed subspace of $L^{2}(\Omega, {\cal F}, P; \C)$  generated by $\{Y_{\pi}, \pi \in \G\}$. Consider the linear mapping $V$ from the complex linear span of $\{Y_{\pi}, \pi \in \G\}$ into $L^{2}_{c}(G, \mu_{Y}):= L^{2}_{c}(G, {\cal B}(G), \mu_{Y}:\C)$ (where the subscript $c$, indicates the restriction to central functions) given by
            $$ V\left(\sum_{j=1}^{n}\alpha_{j}Y_{\pi_{j}}\right) = \sum_{j=1}^{n} \alpha_{j}\chi_{\pi_{j}}.$$
It is straightforward to check that $V$ is isometric. Since the set of all finite linear combinations of characters is dense in $L^{2}_{c}(G, \mu_{Y})$, it follows that $V$ extends to a unitary isomorphism between ${\cal M}$ and $L^{2}_{c}(G, \mu_{Y})$. For each $A \in {\cal B}(G)$, define
$$ \Gamma_{Y}(A) = V^{*}{\bf 1}_{A}.$$
Then it is straightforward to check that $\Gamma_{Y}$ has the desired properties. Moreover, for all $f \in L^{2}(G, \mu_{Y})$,
$$  V\left(\int_{G}f(g)\Gamma_{Y}(dg)\right) = f.$$
Then for all $\pi \in \G$,
$$  V\left(\int_{G}\chi_{\pi}(g)\Gamma_{Y}(dg)\right) =  \chi_{\pi}(g) = VY_{\pi},$$
and we thus obtain (\ref{cram1}). $\hfill \Box$

\vspace{5pt}

We also have a reconstruction theorem of Kolmogorov type:

\begin{theorem} \label{Kol} Given a positive definite function $\Phi: \G \rightarrow \C$ with $\Phi(\epsilon) = 1$, there exists a stationary random field $(Y_{\pi}, \pi \in \GF)$ having covariance $\Phi$.
\end{theorem}

{\it Proof.} By the Bochner theorem of \cite{Hey}, there exists a finite Radon central measure $\mu$ on $(G, {\cal B}(G))$ so that for all $\pi \in \G$,
$$ \Phi(\pi) = \int_{G}\chi_{\pi}(g)\mu(dg),$$
and the normalisation  $\Phi(\epsilon) = 1$ ensures that $\mu$ is a probability measure. Now define $(Y_{\pi}, \pi \in \G)$ on the probability space $(G, {\cal B}(G), \mu)$ by the prescription $Y_{\pi} = \chi_{\pi}$, for each $\pi \in \GF$. Then the field is automatically decomposable and is stationary since
\bean \E(Y_{\pi_{1}}\overline{Y_{\pi_{2}}}) & = & \int_{G}\chi_{\pi_{1}}(g)\overline{\chi_{\pi_{2}}}
(g)\mu(dg)\\
& = & \int_{G}\chi_{\pi_{1} \otimes \pi_{2}^{*}}(g)\mu(dg)\\
& = & \E(Y_{\pi_{1} \otimes \pi_{2}^{*}}\overline{Y_{\epsilon}}), \eean
since $Y_{\epsilon} = \chi_{\epsilon} = 1$. $\hfill \Box$

\vspace{5pt}

Let $\rho$ be a central probability measure on $G$. Then $\widehat{\rho}: \G \rightarrow \C$ is positive definite and so is the covariance of a stationary random field by Theorem \ref{Kol}. So Theorem \ref{Kol} tells us that there are a rich variety of stationary random fields on $\G$. For example, suppose that $G$ is a compact, connected Lie group and that $\rho$ is Gaussian, so that $\rho(d\sigma) = k_{1}(\sigma)d\sigma$, where $(k_{t}, t \geq 0)$ is the heat kernel on $G$. Then for each $\pi \in \G,$
$$ \widehat{\rho}(\pi) := \int_{G}\chi_{\pi}(g)\rho(dg) =  d_{\pi}e^{-\kappa_{\pi}},$$
where $d_{\pi}$ is the dimension of the complex linear space in which $\pi$ acts, and $\{\kappa_{\pi}, \pi \in \G\}$ is the Casimir spectrum. A large class of non-Gaussian infinitely divisible central measures on $G$ may be obtained by subordination of the heat kernel (see e.g. \cite{App3} or Chapter 4 of \cite{App}).

\section{Examples - Time Series}

It is interesting to seek examples in the case where $G$ is a rank-one, connected, compact Lie group. Then the lattice of weights is a subset of the real line, and so inherits an ordering that can be used to develop a theory of time series, by analogy with the familiar one on the group of integers. As an example, let us consider the group $G = SU(2)$. In this case $\G$ is in one-to-one correspondence with the set $\Z_{+}$ (with $0$ corresponding to $\epsilon$) and we may consider the AR$(1)$ process defined for each $n \in \Z_{+}, \lambda \in \C$ by
\begin{equation} \label{AR1}
 Y_{n} = \lambda Y_{n-1} + Z_{n},
\end{equation}
where $Y_{-1}:=0$.

Here we may take the index $n \in \Z_{+}$ as labelling the unique equivalence class of irreducible representations having representation space with dimension $n+1$.

We show that this process cannot be stationary.

Define the backwards shift operator $B$ on the linear space generated by $\{Y_{n}, n \in \Z_{+}\}$ by
$$ BY_{n} = Y_{n-1}. $$
Then
\begin{eqnarray} \label{ARep}
Y_{n} & =  & (I - \lambda B)^{-1}Z_{n} \nonumber \\
& = & \sum_{k=0}^{n}\lambda^{k}Z_{n-k}.
\end{eqnarray}

Note that in contrast to the familiar case of $G = \mathbb{T}$, no condition is needed on $\lambda$ to obtain the moving average representation (\ref{ARep}) as this series is finite.
It follows easily from (\ref{ARep}) that $(Y_{n}, n \in \Z_{+})$ has covariance
\bean \E(Y_{n+h}\overline{Y_{n}}) & = & \sum_{k=0}^{n+h}\sum_{l=0}^{n}\lambda^{k}\overline{\lambda^{l}}\E(Z_{n+h-k}\overline{Z_{n-l}})\\
& = & \overline{\lambda^{-h}}\sum_{k=0}^{n}|\lambda|^{2(k+h)}\\
& = & \left\{ \begin{array}{c c} \lambda^{h}\left(\ds\frac{1-|\lambda|^{2n+ 2}}{1 - |\lambda|^{2}}\right)&~\mbox{if}~|\lambda| \neq 1 \\
 & \\
(n+1)\overline{\lambda^{-h}}&~\mbox{if}~|\lambda| = 1. \end{array} \right.
\eean

\vspace{5pt}

On the other hand, consider the MA$(q)$ process on $\widehat{SU(2)}$ given by
$$ Y_{n} = \sum_{k=0}^{q}\beta_{k}Z_{n-k},$$

where $\beta_{k} \in \C$ for all $k \in \Z_{+}$. Then by standard arguments, $(Y_{n}, n \in \Z_{+})$ is easily seen to be a stationary random field with covariance function

$$ \E(Y_{n+h}\overline{Y_{n}}) = \left\{ \begin{array}{c c} \sum_{k=0}^{q-h}\beta_{k+h}\overline{\beta_{k}} &~\mbox{if}~0 \leq h \leq q \\
& \\
0  &~\mbox{if}~h > q.  \end{array} \right. $$

\section{The Hypergroup Connection}

Let $K$ be a non-empty locally compact Hausdorff space which is equipped with an involution $x \rightarrow x^{\prime}$. Let $M^{b}(K)$ be the complex linear space of all bounded, complex Radon measures on $K$. We say that $(K, *)$ is a {\it hypergroup} if there is a binary operation $*$ defined on $M^{b}(K)$ with respect to which $M^{b}(K)$ is a algebra, and if certain axioms hold. We state only one of these here; that there must exist a ``neutral element'' $e \in K$ so that for all $x \in K$,
$$ \delta_{x} * \delta_{e} = \delta_{e} * \delta_{x} = \delta_{x},$$
where $\delta_{x}$ is the Dirac mass at $x$. The others may be found on p.9 of \cite{BH}; they will play no direct role in the sequel. Examples are locally compact groups (where $*$ is the usual convolution of measures), double coset spaces and the unitary dual of a compact group (see below).
The hypergroup is said to be {\it discrete} if $K$ is equipped with the discrete topology, and {\it commutative} if
$ \mu_{1} * \mu_{2} = \mu_{2} * \mu_{1}$ for all $\mu_{1}, \mu_{2} \in M^{b}(K)$.

Now let $\G$ be the unitary dual of a compact group $G$. It becomes a discrete, commutative hypergroup, with neutral element $\epsilon$ and involution $\pi \rightarrow \pi^{*}$, under the convolution:

\begin{equation} \label{conv}
\delta_{\pi_{1}} * \delta_{\pi_{2}} = \ds\sum_{\pi \in \G}M(\pi_{1} \otimes \pi_{2}, \pi)\delta_{\pi},
\end{equation}

for each $\pi_{1}, \pi_{2} \in \G,$ relative to the decomposition $$\pi_{1} \otimes \pi_{2} = \bigoplus_{\pi \in \G}M(\pi_{1} \otimes \pi_{2}, \pi)\pi.$$ The convolution (\ref{conv}) is extended to general measures in $M^{b}(\G)$ by taking weak limits of linear combinations.  Note that this is not the same convolution as that given in \cite{BH}, p.13, where the following is found:

\begin{equation} \label{conv1}
\delta_{\pi_{1}} *^{\prime} \delta_{\pi_{2}} = \ds\sum_{\pi \in \G}\frac{d_{\pi}}{d_{\pi_{1}}d_{\pi_{2}}}M(\pi_{1} \otimes \pi_{2}, \pi)\delta_{\pi},
\end{equation}
with $d_{\pi}$ being the dimension of the representation space corresponding to $\pi \in \G$.

Following section 8.2 in \cite{BH}, the survey article \cite{Hesurv} and the original source \cite{LaLe} we define a {\it stationary random field} over a commutative hypergroup $K$ to be a mapping $X:K \rightarrow L^{2}(\Omega, {\cal F}, P; \C)$ which has covariance
$$ C(a,b) = \E(X_{a}\overline{X_{b}}),$$
that satisfies the stationarity condition:
\begin{equation} \label{stathyp}
C(a,b) = \int_{K}C(x, e)(\delta_{a} * \delta_{b^{\prime}})(dx),
\end{equation}
for each $a,b \in K$.

Now suppose that $(Y_{\pi}, \pi \in \G)$ is a stationary random field on $\G$ in the sense of (\ref{statdef}). We will show that it is also stationary in the hypergroup sense, by using the convolution (\ref{conv}) to define the hypergroup structure. Note that if we used (\ref{conv1}) then this assertion would be false. It is enough to show that (\ref{stathyp}) is satisfied. Indeed for all $\pi_{1}, \pi_{2} \in \G$,
\bean C(\pi_{1}, \pi_{2}) & = & \E(Y_{\pi_{1}}\overline{Y_{\pi_{2}}})\\
& = & \E(Y_{\pi_{1} \otimes \pi_{2}^{*}}\overline{Y_{\epsilon}})\\
& = &C(\pi_{1} \otimes \pi_{2}^{*}, \epsilon) \\
& = & \sum_{\pi \in \G}M(\pi_{1} \otimes \pi_{2}^{*}, \pi)C(\pi, \epsilon)\\
& = & \int_{\G}C(\pi, \epsilon)(\delta_{\pi_{1}} * \delta_{\pi_{2}^{*}})(d\pi),
\eean
as required.

\vspace{5pt}

Let ${\cal M}_{Y}(\G)$ be the closure in $L^{2}(\Omega, {\cal F}, P; \C)$ of $\{Y_{\pi}, \pi \in \G\}$. Following \cite{LaLe, Leit, Hesurv}, for fixed $\pi^{\prime} \in \G$, we define the {\it translation operator $\tau_{\pi^{\prime}}$ associated to $Y$} to be the linear contraction in ${\cal M}_{Y}(\G)$ obtained by continuous linear extension of the prescription

\begin{equation} \label{trans}
\tau_{\pi^{\prime}}(Y_{\pi}) = Y_{\pi \otimes \pi^{\prime}},
\end{equation}

\noindent for each $\pi \in \G$. Note that because of the rather concrete context in which we work, the definition (\ref{trans}) is much more transparent than that in the general hypergroup case. A useful list of properties of such operators is collected in Theorem 2 of Leitner \cite{Leit}.

Now let ${\cal A}$ be a family of subsets of $\G$. For each $A \in {\cal A}$, let ${\cal M}_{Y}(A)$ be the closure of the linear span of $\{Y_{\pi}, \pi \in A\}$, and ${\cal M}_{Y}: = \bigcap_{A \in {\cal A}}{\cal M}_{Y}(A)$. We say that the stationary field $Y$ is {\it ${\cal A}$-singular} if
${\cal M}_{Y} = {\cal M}_{Y}(\G)$, {\it ${\cal A}$-regular} if ${\cal M}_{Y} = \{0\}$, and {\it ${\cal A}$-adapted} if $\tau_{\pi}({\cal M}_{Y}) \subseteq {\cal M}_{Y}$ for all $\pi \in \G$. The following abstract version of the {\it Wold decomposition} is proved for general commutative hypergroups in Theorem 2.2.5.2 of \cite{Hesurv}. We will be content to state the result.

\begin{theorem}[The Wold Decomposition] \label{HeyWold}
If $Y$ is an ${\cal A}$-adapted stationary random field, then there is a unique orthogonal decomposition
\begin{equation}  \label{HeyWold1}
Y_{\pi} = Y^{(1)}_{\pi} + Y^{(2)}_{\pi},
\end{equation}
for all $\pi \in \G$, where $Y^{(1)}$ is ${\cal A}$-regular, and $Y^{(2)}$ is ${\cal A}$-singular.
\end{theorem}

\vspace{5pt}

If $G$ is a rank one, connected, compact Lie group then it is natural to choose ${\cal A}$ in accordance with the lattice structure, e.g. for $G = SU(2), {\cal A} = \{A_{n}, n \in \Z_{+}\}$, where $A_{n}: = \{0, 1, \ldots, n\}$.

\vspace{5pt}

We conjecture that there is a generalisation to this context of the classical result that can be found e.g. in Chapter 4 of \cite{Lamp}, whereby the the absolutely continuous measure $\mu_{1}$ and the singular measure $\mu_{2}$ which arise in the Lebesgue decomposition of the spectral measure of $Y$ are themselves the spectral measures of the processes $Y^{(1)}$ and $Y^{(2)}$ (respectively) of (\ref{HeyWold1}).

\vspace{5pt}

{\bf Acknowledgement}. The author would like to thank N.H.Bingham for stimulating discussions about prediction, in relation to the ideas in \cite{Bing}, that inspired the work of this note. He is also indebted to H.Heyer for very helpful comments.

\end{document}